\documentclass[12pt,letter]{amsart}

\title[Monodromies of a degeneration of irreducible symplectic manifold]{On monodromies of a degeneration of irreducible symplectic K\"ahler manifolds}
\author{Yasunari Nagai}

\address{ 
Korea Institute for Advanced Study (KIAS), 
207-43 Cheongnyangni 2-dong, Dongdaemun-gu, Seoul 130-722, Korea } 
\email{nagai@kias.re.kr}

\usepackage{amsthm}
\usepackage{amssymb}
\usepackage{amsxtra}
\usepackage[mathscr]{eucal}
\usepackage{enumerate}
\usepackage[all]{xy}

\usepackage{times,mathptmx}

\theoremstyle{plain}
\newtheorem{theorem}{Theorem}[section]
\newtheorem{lemma}[theorem]{Lemma}
\newtheorem{proposition}[theorem]{Proposition}
\newtheorem{corollary}[theorem]{Corollary}
\newtheorem{conjecture}[theorem]{Conjecture}
\newtheorem*{theorem*}{Theorem}
\newtheorem*{corollary*}{Corollary}
\theoremstyle{definition}
\newtheorem{definition}[theorem]{Definition}
\newtheorem{example}[theorem]{Example}

\theoremstyle{remark}
\newtheorem{remark}[theorem]{Remark}
\newtheorem*{acknowledgement}{Acknowledgemenet}

\def\lto{\longrightarrow}

\DeclareMathOperator{\supp}{supp}

\DeclareMathOperator{\Hilb}{Hilb}
\DeclareMathOperator{\Sym}{Sym}
\DeclareMathOperator{\Kum}{Kum}

\DeclareMathOperator{\id}{id}

\DeclareMathOperator{\nilp}{nilp}
\DeclareMathOperator{\Gr}{Gr}
\DeclareMathOperator{\Ker}{Ker}
\DeclareMathOperator{\iM}{Im}
\DeclareMathOperator{\rank}{rank}
\DeclareMathOperator{\Sing}{Sing}
\DeclareMathOperator{\rE}{Re}

\begin{document}

\baselineskip 15pt
\parskip 5pt

\maketitle

\begin{abstract}
 We study the monodromy operators on the Betti cohomologies associated to
 a good degeneration of irreducible symplectic manifold and we show that
 the unipotency of the monodromy operator on the middle cohomology is at
 least the half of the dimension. This implies that the ``mildest''
 singular fiber of a good degeneration with non-trivial monodromy 
 of irreducible symplectic manifolds is quite different 
 from the generic degeneration of abelian varieties 
 or Calabi-Yau manifolds.
\end{abstract}


\section*{Introduction}

For the study of smooth algebraic (or analytic) varieties, it is often
important and useful to consider its \emph{degenerations}. 
One can see this principle, for example, by remembering the role played by 
the singular fibers in the theory of elliptic surfaces due to Kodaira.

The Kodaira singular fibers are completely described in terms of the
periods and monodromies around the singular fiber. 
The study of the periods and degenerations of abelian varieties is one
of the most direct generalizations of 
the theory of elliptic surfaces to higher dimensions. 
The theory of periods is generalized by Griffiths and many
other contributers to the general situation using the \emph{variation of
Hodge structures}. Since then, many significant results are proved using
this mechanism. In some cases, the Hodge structure on the middle
cohomology, i.e. $H^n(X,\mathbb C)$ for $X$ of dimension $n$, plays an
important role. As an example, we can recall 
the classification of Kulikov models for
degeneration of K3 surfaces \cite{Ku,P-P} and the proof of global
Torelli theorem for K3 surfaces via Kulikov models by Friedman
\cite{Fr}. 

An \emph{irreducible symplectic K\"ahler manifold} is a generalization of a K3
surface to higher dimension. As is well known, the local Torelli
theorem holds for irreducible symplectic manifolds on the
\emph{second} cohomology group (not only on the middle cohomology). 
From this result, one can easily suppose some strong similarities between 
the irreducible symplectic manifolds and K3 surfaces. 
There is no reason that prevents us from studying degenerations of
irreducible symplectic manifolds. 

But if we consider the period and 
monodromy only on the second cohomology as invariants of the
degeneration of irreducible symplectic manifolds, it seems that one
misses some important information. Since the weight of
the Hodge structure on the second cohomology is 2, the index of
unipotency of the monodromy is either 0, 1, or 2. 
The monodromy operator should have some information about the
``complexity'' of the degenerate fiber. In the case of Kulikov models of
K3 surfaces or toroidal degenerations of abelian varieties, the
unipotency of the monodromy on the middle cohomology certainly
corresponds to the combinatorial complexity of the singular fiber. 
There are $(n+1)$ patterns of local combinations of 
the components of the degenerate variety with normal crossings of
dimension $n$, so the unipotency of the monodromy on the second
cohomology, \emph{a priori}, may not fully reflect 
the combinatorial information of the degenerate fiber. 

In this article, we examine the relation between 
the monodromy operators on the cohomologies, in particular the relation
between the monodromies on the second cohomology and the middle
cohomology. 

To consider this problem in a general situation, we assume
a semi-stable model which is called a \emph{good degeneration} (see
Definition \ref{gooddegsympl}). 
In fact, we can construct an example of a good degeneration 
(Theorem \ref{exampleofdsv}). 
Our main result is the following.

\begin{theorem*}[Theorem \ref{main}]
 Let $\pi:\mathcal X\to \Delta$ be a good degeneration of irreducible 
 symplectic $2n$-folds. Let $T_{2n}$ be the monodromy operator on 
 $H^{2n}(\mathcal X_t,\mathbb C)\;(t\neq 0)$ associated to the
 family $\pi$, and $N_{2n}=\log T_{2n}$. Assume $N_{2n}^n=0$, then
 $N_{2n}=0$.
\end{theorem*}

The maximal $l$ with $N_{2n}^l\neq 0$ is called the unipotency of
$T_{2n}$. This theorem asserts that the unipotency of the monodromy on the
middle cohomology associated to a good degeneration of irreducible
symplectic manifold is 0 or not less than $n$. This is quite different
from the situation one can expect in the case of 
general semi-stable degeneration. 
For example, we can easily show the following corollary. 

\begin{corollary*}[See Corollary \ref{nocycle}]
 Consider a good degeneration of irreducible
 symplectic $2n$-folds with non-trivial monodromy on the middle cohomology, 
 and let $\Gamma$ be the dual graph of the configuration of the irreducible
 components of the singular fiber.
 Then the dimension of the topological realization $|\Gamma |$ is at least $n$.
\end{corollary*}

\noindent
This phenomenon can be seen as an aspect of 
the general principle that the geometry of
irreducible symplectic manifold is very restrictive. 

\paragraph{{\it Plan of the article}}
This article consists of five sections: 
In the first section, we review some
necessary definitions about symplectic K\"ahler manifolds, degenerations
and monodromy associated to them. In \S 2, we consider an example of
degeneration of irreducible symplectic manifold 
arising from a family of K3 surfaces by the Hilbert scheme
construction and compute the monodromies associated to it. 
In fact, this example motivates this research.
The next section is devoted to the study of a family of 
generalized Kummer varieties and its associated monodromy. 
In \S 4, we define a notion of a good degeneration of
symplectic K\"ahler manifold. 
In the last section, we state a conjecture based on the 
examples in \S\S 2,3, and get a partial answer to the 
conjecture for good degenerations of irreducible symplectic 
manifold. 

\begin{acknowledgement}
 The author would like to thank Prof. Eyal Markman for pointing out a 
 mistake in \S 2 of the earlier version and Dr. Jaeyoo
 Choy for his comments.
 He is also grateful to the referee for his careful reading and 
 appropriate comments. The author is financially 
 supported by PD Fellowship of
 Korea Institute for Advanced Study (KIAS). 
\end{acknowledgement}

\section{Definitions and notations}

In this section, 
we collect some basic definitions and notations which are
needed in this article.

\subsection*{Symplectic K\"ahler manifold}

Let us start with a review of some basics of compact symplectic
K\"ahler manifold. A fundamental reference is \cite{Be}. See also
\cite{Huy1,Huy2}. 

\begin{definition}\label{IS}
 Let $X$ be a compact K\"ahler manifold of dimension $2n$.
 A holomorphic $2$-form $\sigma \in H^0(X,\Omega _X^2)$ 
 on $X$ is a (holomorphic) \emph{symplectic form} if its top exterior 
 power $\sigma ^{\wedge n}\in H^0(X,\Omega _X^{2n})$ 
 is nowhere vanishing. The pair $(X,\sigma _X)$ is called a
 holomorphic \emph{symplectic K\"ahler manifold}. 
 $X$ is said to be 
 an \emph{irreducible symplectic} if the following conditions are
 satisfied. 
 \begin{enumerate}[(i)]
  \item There exists a holomorphic symplectic form $\sigma$.
  \item The space of global holomorphic 2-forms $H^0(X,\Omega _X^2)$ is
	spanned by the symplectic form $\sigma$. 
  \item $X$ is simply connected.
 \end{enumerate}
 
 Given an irreducible symplectic manifold $X$, we have a 
 non-degenerate primitive quadratic form $q_X$ on $H^2(X,\mathbb Z)$, 
 which is called the \emph{Beauville-Bogomolov form}. 
\end{definition}

Famous decomposition theorem for compact K\"ahler manifolds with trivial
canonical bundle implies that a compact symplectic K\"ahler manifold is
an \'etale quotient of a product of a complex torus and irreducible
symplectic manifolds. 
There are few known examples of irreducible symplectic K\"ahler
manifold. Here we describe two families of examples, 
which are classically known.

\begin{example}[\cite{Be}]
 Let $S$ be a K3 surface and $\Hilb ^n(S)$ the Hilbert scheme
 (or Douady space) of $0$-dimensional sub-schemes of length $n$ on $S$. 
 Then, $\Hilb ^n(S)$ is an irreducible symplectic K\"ahler 
 manifold of dimension $2n$. 
 For the second Betti cohomology of $\Hilb ^n(S)$, we have an
 isomorphism 
 \begin{equation}\label{BBhilb}
  H^2(\Hilb ^n(S),\mathbb C)\cong H^2(S,\mathbb C)\oplus \mathbb C\cdot
 [E], 
 \end{equation}
 where $E$ is an irreducible exceptional divisor of the birational 
 morphism 
 \[
 \Hilb ^n(S)\to \Sym ^n(S), 
 \]
 which is called the Hilbert--Chow morphism (\cite{Be}) and decomposition
 in \eqref{BBhilb} is orthogonal with respect to the
 Beauville--Bogomolov form.
 
 More generally, a compact connected component $M$ of the moduli 
 space of stable sheaves on a K3 surface $S$ is known to be 
 an irreducible symplectic manifold which is deformation equivalent to
 $\Hilb ^n(S)$ of appropriate dimension. An irreducible symplectic
 manifold $X$ is said to be \emph{of Hilbert type} if $X$ 
 is deformation equivalent to some $\Hilb ^n(S)$.
\end{example}

\begin{example}[\cite{Be}]\label{kummer}
 Let $A$ be a complex torus and consider $\Hilb ^{n+1}(A)$ and its
 Albanese morphism $\alpha :\Hilb ^{n+1}(A)\to A$, which is a locally
 trivial fiber bundle. Take its fiber $\Kum ^n(A)=\alpha ^{-1}(0)$. Then
 this is an irreducible symplectic K\"ahler manifold of dimension $2n$. 
 We can regard $\Kum ^n(A)$ as a resolution of $A^n/\mathfrak
 S_{n+1}$. Let $E$ be the exceptional divisor. Then $E$ is irreducible
 and we have 
 \[
 H^2(\Kum ^n(A),\mathbb C)\cong H^2(A,\mathbb C)\oplus \mathbb C\cdot [E]
 \]
 and this is also orthogonal with respect to the Beauville-Bogomolov form.
\end{example}

More recently, O'Grady
constructed two sporadic examples of irreducible symplectic manifold
\cite{OG1,OG2}. It is noteworthy that these four types of examples 
are all of the known examples for the moment.

\subsection*{Degeneration and Monodromy}

\begin{definition}
 Let $Y$ be a compact K\"ahler manifold and $\Delta$ a unit disk in $\mathbb
 C$. A flat proper morphism $\pi:\mathcal X\to \Delta$ of a normal
 complex K\"ahler space $\mathcal X$ is said to be a \emph{degeneration} 
 or \emph{degenerating family} of $Y$ 
 if $\pi$ is smooth over $\Delta-\{0\}$
 and $\mathcal X_t=\pi ^{-1}(t)\cong Y$ for some $t\in \Delta$.
 The fiber $X=\mathcal X_0=\pi ^{-1}(0)$ is said to be the \emph{singular
 fiber} if $\pi$ is not smooth.
\end{definition}

Given a degenerating family $\mathcal X$, 
we have the monodromy operators on the cohomologies of $\mathcal X_t\;
(t\neq 0)$.

\begin{definition}
 Let $\pi :\mathcal X\to \Delta$ be a degenerating family. 
 Then parallel displacement on the local system 
 $R^m(\pi _0)_* \mathbb C$ on $\Delta-\{0\}$ 
 induces a homomorphism
 \[
  T_m:\pi _1(\Delta-\{0\})\to GL(H^m(\mathcal X_t,\mathbb C)), 
 \]
 the \emph{monodromy representation}. We denote also by $T_m$ 
 the image of a generator under the monodromy representation, which is called
 the \emph{monodromy operator}.
\end{definition}

By the monodromy theorem, $T_m$ is \emph{quasi-unipotent}, i.e. 
$(T_m^k-I)^N=0$ for some $k,N\in \mathbb N$. 
Therefore, changing the base by a
cyclic cover $\Delta \to \Delta,\; t\mapsto t^k$, we can always 
make the monodromy operator $T_m$ \emph{unipotent}, i.e., $k=1$.  

\begin{definition}
 Let $T$ be an unipotent automorphism of a finite dimensional vector
 space over a field of characteristic zero.
 The logarithm of $T$ is defined by 
 \[
  N=\log T=(T-I)-(T-I)^2+\cdots +(-1)^{n+1}(T-I)^n +\cdots . 
 \]
 Note that the right hand side is a finite sum because $T$ is
 unipotent and the logarithm $N$ is a nilpotent endomorphism. 
 Of course, we can reconstruct $T$ from $N$ by the exponential:
 \[
 T=\exp N=I+\frac{N}{1!}+\frac{N^2}{2!}+\cdots +\frac{N^n}{n!}+\cdots .
 \]
 We define the \emph{index of nilpotency} of $N$ by
 \[
 \nilp (N)=\max \{k \mid N^k\neq 0\}. 
 \]
 We mean by the \emph{index of unipotency} 
 of a unipotent $T$ the index of nilpotency of $N=\log T$.  
\end{definition}

If $T_m$ is the unipotent monodromy operator on 
$H^m(\mathcal X_t,\mathbb C)$
associated to a degeneration $\pi:\mathcal X\to \Delta$, 
then by the theory of variation of Hodge structures, 
more precisely, by $SL_2$-orbit theorem \cite{Sc} for example,  
we know that the index of unipotency of $T_m$ is at most $m$, i.e.
\begin{equation}\label{boundofnilp}
 0\leq \nilp (N_m)\leq m
\end{equation}
for $N_m=\log T_m$ (see also a lecture note by Griffiths \cite{Topics},
Chapter IV).

\section{Example: the case of Hilbert type}

In this section, we consider degenerations of irreducible symplectic 
manifolds of Hilbert type, i.e. irreducible symplectic manifolds which
is deformation equivalent to $\Hilb ^n(S)$ for some K3 surface $S$.

Let us consider an easy example.

\begin{example}\label{HilbFam}
 Let $p : \mathcal S\to \Delta $ be a projective degeneration of 
 K3 surfaces. Consider the Hilbert scheme$\overline{\mathcal X}_n
 = \Hilb ^n(\mathcal S/\Delta)$ of 
 $0$-dimensional sub-schemes relative to $p$, 
 and take the normalization $\mathcal X_n \to \overline{\mathcal X}_n$. 
 Then, the natural morphism $\pi_n : \mathcal X_n\to \Delta$ 
 is projective and $\mathcal X_n$ is a normal K\"ahler space after
 shrinking $\Delta$ if necessary. Therefore, $\pi_n$  
 is a degeneration of irreducible symplectic manifolds, 
 whose general fiber $(\mathcal X_n)_t$ is isomorphic to 
 $\Hilb ^n(\mathcal S_t)$. 
\end{example}

Consider the monodromy operator on the cohomologies associated to this
degeneration. 
Let us assume that the monodromy operator 
$T'$ on the second cohomology group $H^2(\mathcal S_t,\mathbb
C)\; (t\neq 0)$ associated to the family $p:\mathcal S\to \Delta$ is
unipotent and let $N'=\log T'$. Then by
\eqref{boundofnilp}, we have 
\[
 \nilp (N')=0,\, 1,\mbox{ or } 2. 
\]
Since we have 
$H^2(\Hilb ^n(\mathcal S_t),\mathbb C)=H^2(\mathcal S_t,\mathbb
C)\oplus \mathbb C\cdot [E_t]$ on the general fiber and the class 
$[E_t]$ is clearly invariant under the action of the monodromy $T_2$ 
on the second cohomology $H^2(\Hilb ^n(\mathcal S_t),\mathbb C)$, 
we have 
\[
 T_2=T'\oplus \id, 
\]
in particular, $T_2$ is unipotent. In fact we can say more about the
monodromy operator $T_m$ on 
$H^m((\mathcal X_n)_t,\mathbb C)$:

\begin{proposition}\label{MHilbFam}
 Notation as above. Then for $m\leqslant 2n$, 
 \begin{enumerate}[(i)]
  \item $T_m$ is unipotent. 
  \item Let $N_{2k}=\log T_{2k}$. Then, 
	$\nilp (N_{2k})=k\cdot \nilp (N_2)$ for
	$k\leqslant n$. In particular 
	$\nilp (N_{2k})\in \{0,k,2k\}$. 
\end{enumerate}
\end{proposition}

Let $SH^*(\Hilb ^n(S),\mathbb C)$ be the sub-algebra of 
$H^*(\Hilb ^n(S),\mathbb C)$
generated by $H^2(\Hilb ^n(S),\mathbb C)$. For 
$SH^*(\Hilb ^n(S),\mathbb C)$, we have the following general 
result.

\begin{proposition}[Verbitsky, \cite{Bo}. See also \cite{Huy2}
 \S 24]\label{Verb}
 Let $X$ be an irreducible symplectic manifold of dimension $2n$ and 
 let $SH^*(X,\mathbb C)$ be the sub-algebra of $H^*(X,\mathbb C)$
 generated by $H^2(X,\mathbb C)$. Then 
 \[
 SH^*(X,\mathbb C)\cong 
 \Sym ^* H^2(X,\mathbb C) /\langle \alpha ^{n+1}\mid
 q_X(\alpha)=0 \rangle
 \]
 where $q_X$ is the Beauville--Bogomolov form (Definition \ref{IS}). 
\end{proposition}

By this proposition, we have a natural injection 
\[
 \Sym ^k H^2(\Hilb ^n(S),\mathbb C)\hookrightarrow 
 H^{2k}(\Hilb ^n(S),\mathbb C)
\] 
for $k\leqslant n$. 
The easiest case to prove the proposition is the case where $n=2$, 
i.e., the case of $\Hilb ^2(S)$. In this case, one can easily see that 
$H^*(\Hilb ^2(S))=SH^*(\Hilb^2(S))$. This implies that
$H^4(\Hilb ^2(S))=\Sym ^2H^2(\Hilb ^2(S))$. Therefore, the proposition
is just a consequence of following lemma.

\begin{lemma}\label{tensornilp}
 Let $V_1$, $V_2$ be finite dimensional vector spaces 
 over a field of characteristic zero, $T_i$ be a
 unipotent automorphism on $V_i$, and $N_i=\log T_i$. Then $T_1\otimes
 T_2$ is also unipotent and 
 \[
 \log (T_1\otimes T_2) = N_1\otimes I + I \otimes N_2
 \]
 on $V_1\otimes V_2$. Moreover, we have  
 \[
  \nilp (\log (T_1\otimes T_2))=\nilp(N_1)+\nilp (N_2).
 \]
In particular, $\Sym ^k T_1$ on $\Sym ^k V_1$ is unipotent and 
 \[
 \nilp (\log (\Sym ^k T_1))=k\cdot \nilp (N_1).
 \]
\end{lemma}

\begin{proof}
 The first assertion is just a property of exponentials and logarithms
 of the matrices. The second equality follows from 
 \[
 (N_1\otimes I + I\otimes N_2)^k(v_1\otimes v_2)
 =\sum _{i=0}^k \begin{pmatrix} k \\ i \end{pmatrix} N_1^k(v_1)\otimes
 N_2^{k-i}(v_2). 
 \]
 Since $\Sym ^k V_1\subset V_1^{\otimes k}$, we have 
 \[
 \nilp (\log (\Sym ^k T_1))\leqslant 
 \nilp (\log (T_1^{\otimes k}))=k\cdot \nilp (N_1)
 \] 
 On the other hand, there exists $v_1\in V_1$ such that $N_1^l(v_1)\neq 0,\;
 N_1^{l+1}(v_1)=0$ for $l=\nilp (N_1)$. Therefore
 \[
 (\log (\Sym ^k T_1))^{kl}(v_1^k)=M_{k,l}\cdot (N^l(v_1))^k\neq 0 
 \]
 for some positive integer $M_{k,l}$, 
 which shows the last claim.
\end{proof}

Since $H^{2k}(\Hilb ^n(S))$ is not generated by
$H^2(\Hilb ^n(S))$ for $n>2$, we need some more argument to prove Proposition 
\ref{MHilbFam} for $n>2$. 
This can be done by applying G\"ottsche--Soergel formula.

\begin{theorem}[G\"ottsche--Soergel \cite{GS}, 
 Theorem 2]\label{GSformula0} 
 Let $S$ be a smooth algebraic surface and 
 fix a natural number $n$. Consider the set $P(n)$ of partitions of
 $n$, i.e., 
 \[
  P(n)=\{\alpha=(\alpha _1, \alpha _2, \cdots ,\alpha _n)
 \mid \alpha _1\cdot 1+\alpha _2\cdot 2 +\cdots +\alpha _n\cdot n=n\},
 \] 
 and put $|\alpha |=\sum _i \alpha _i$. 
 Then, there is a canonical isomorphism
 \begin{equation}\label{GS}
 H^{i+2n}(\Hilb ^n(S))=\bigoplus _{\alpha\in P(n)} 
 \left( H^{i+2|\alpha |} (S^{(\alpha)}) \right),
 \end{equation}
 where $S^{(a)}=\Sym ^a(S)$ for a positive integer $a$ and 
 $S^{(\alpha)}=S^{(\alpha _1)}\times \cdots \times S^{(\alpha _n)}$
 for $\alpha =(\alpha _1,\cdots ,\alpha _n)$ (we regard $S^{(0)}=
 \Sym ^0(S)$ as a one point set). 
\end{theorem}

\begin{remark}
 By this theorem, we know that $H^m(X,\mathbb C)=0$ for odd $m$ if 
 $X$ is deformation equivalent to some $\Hilb ^n(\mbox{K3})$. 
 Therefore, we have nothing to do with the monodromy operators 
 $T_m$ for odd $m$ in this section. 
\end{remark}

\begin{proof}[Proof of Proposition \ref{MHilbFam}]
 Apply Theorem \ref{GSformula0} to the case $S$ is a projective K3
 surface. Since $H^1(S)=H^3(S)=0$, the theorem and K\"unneth formula 
 implies that the cohomology group $H^{2k}(\Hilb ^n(S))$ is a direct sum 
 of tensor products of several symmetric products of $H^2(S)$. 
 G\"ottsche--Soergel decomposition is invariant under the monodromy
 operator of the family $\Hilb ^n(\mathcal S/\Delta)\to \Delta$ 
 because the decomposition is induced by the action of the
 symmetric group on $S^n$. This proves (i). 
 Since the weights of the components of the
 G\"ottsche--Soergel decomposition of $H^{2k}(\Hilb ^n(S))$ do not
 exceed $2k$, (ii) follows from Lemma \ref{tensornilp}.
\end{proof}

Under the assumption that the monodromy
operators $T_m$ on $H^m(\Hilb ^n(S), \mathbb C)$ are unipotent, 
one can generalize (ii) in Proposition \ref{MHilbFam} to
the following theorem. 

\begin{theorem}\label{MofHilbFam}
 Let $\pi:\mathcal X\to \Delta$ be a degenerating family 
 of irreducible symplectic
 manifolds such that $\mathcal X_t\cong \Hilb ^n(S)$ for some $t\neq 0$
 and a projective K3 surface $S$. 
 Let $T_m$ be the monodromy operator 
 on $H^m(\mathcal X_t,\mathbb C)$ 
 associated to this family and assume that $T_{2k}$ is unipotent 
 for $k\leqslant n$. Then, $\nilp (N_{2k})=k\cdot \nilp (N_2)$, 
 where $N_{2k}=\log T_{2k}$. 
 In particular $\nilp (N_{2k})\in \{0,k,2k\}$.
\end{theorem}

In fact, the monodromy action on the cohomology ring of an irreducible
symplectic manifolds of Hilbert type is extensively studied by Markman
\cite{M1,M2}. He constructed generators of the cohomology ring using
the universal family and clarified the action of the group of all
possible monodromy operators which need not be unipotent. Markman proved
the existence of the following monodromy invariant decomposition which
is more or less analogous to the G\"ottsche--Soergel decomposition in the
case of Proposition \ref{MHilbFam}:

\begin{theorem}[Markman, \cite{M2}, Corollary 4.6 and Lemma 4.8]
\label{Markman}
\mbox{ }
\begin{enumerate}[(i)]
 \item Let $X$ be an irreducible symplectic manifold, 
       $A_l$ the sub-algebra of $H^*(X,\mathbb C)$ 
       generated by $\oplus _{i=0}^lH^i(X,\mathbb C)$ and
       $[A_l]^j=A_l\cap H^j(X,\mathbb C)$. Then, we have a monodromy
       invariant decomposition
       \[
       H^j(X,\mathbb C)=[A_{j-2}]^j \oplus C_j 
       \]
       for some subspace $C_j\subset H^j(X,\mathbb C)$. 
 \item If $X=\Hilb ^n(S)$, $C_{2k}$ has a monodromy invariant
       decomposition 
       \[
       C_{2k}=C'_{2k}\oplus C''_{2k}
       \] 
       such that $C'_{2k}$ is a one dimensional character and 
       $C''_{2k}$ is $H^2(\Hilb ^n(S))\otimes \chi ''$ where 
       $\chi ''$ is a one dimensional character with values $\{\pm 1\}$
       as representations of monodromy, 
       unless $C'_{2k}$ or $C''_{2k}$ does not vanish. 
\end{enumerate}
\end{theorem}

We should note that the part (i) of this theorem is based on 
the result of Verbitsky and Looijenga--Lunts \cite{V1,V2,LL} about the action of
$\mathfrak{so}(4, b_2-2)$ (where $b_2$ is the second Betti number) 
on the cohomology ring arising from the
existence of the hyper-K\"ahler metric.

Theorem \ref{MofHilbFam} easily follows form this theorem.

\begin{proof}[Proof of Theorem \ref{MofHilbFam}]
 We proceed by induction on $k$. Assume that $\nilp (N_{2l})=l\cdot
 \nilp (N_2)$ for $l<k$. Then we have $\nilp
 (N_{2k|[A_{2k-2}]^{2k}})\leqslant k\cdot \nilp (N_2)$. 
 Moreover, $\Sym ^kH^2(\mathcal X_t, \mathbb C)$ injects to
 $[A_{2k-2}]^{2k}$ (Proposition \ref{Verb}). This implies that $\nilp
 (N_{2k|[A_{2k-2}]^{2k}})= k\cdot \nilp (N_2)$ by Lemma \ref{tensornilp}.
 On the other hand, (ii) of Theorem \ref{Markman} implies that 
 the nilpotency of $N_{2k|C_{2k}}$ is at most $\nilp (N_2)$. 
 Since the decomposition of Theorem \ref{Markman}, (i) is monodromy
 invariant, we conclude that $\nilp (N_{2k})=k\cdot \nilp (N_2)$. 
\end{proof}

\section{Example: A family of generalized Kummer varieties}

From the results of Proposition 
\ref{MHilbFam} and Theorem \ref{MofHilbFam}, 
it is natural to ask whether there is some restriction on the index of
unipotency of the monodromy operator for more general case.
The life in general is not as simple as in the case of Hilbert type, 
since not only the cohomology ring of 
a general irreducible symplectic K\"ahler manifold
is not generated by the second degree part, but also
it can have non-zero odd degree, so some mysterious thing may happen 
in the higher degree cohomologies. 
To catch a glimpse of these general cases, we may use the family
of generalized Kummer varieties as a test case. 

\begin{definition}\label{kumfam}
 Let $p:\mathcal A\to
 \Delta ^*=\Delta - \{0\}$ be a smooth projective family of abelian surfaces 
 with $0$-section and consider its relative Hilbert scheme of $(n+1)$
 points $\tilde {\pi}:\Hilb ^{n+1}(\mathcal A/\Delta^*)\to \Delta ^*$. 
 Then we have the commutative diagram
 \[
 \xymatrix @C=8pt { 
 \Hilb ^{n+1}(\mathcal A/\Delta ^*)\ar[dr]_{\tilde \pi} \ar[rr]^{\quad\alpha} && 
 \mathcal A\ar[dl]^p\\
 & \Delta ^*
 } 
 \]
 where $\alpha$ is the Albanese morphism over $\Delta ^*$, and let 
 $\Kum ^n(\mathcal A/\Delta^*)$ be the inverse image of the $0$-section
 of $p$ by $\alpha$. Then $\pi :\Kum ^n(\mathcal A/\Delta ^*)\to \Delta ^*$ is a
 family of \emph{generalized Kummer varieties} of dimension $2n$ (cf.
 Example \ref{kummer}). 
\end{definition}

\begin{remark}\label{compactify}
 We should be able to compactify $\pi :\Kum ^n(\mathcal A/\Delta^*)\to \Delta^*$ by
 allowing some mild singular fiber of $p$ over the origin as in the case
 of Hilbert schemes in \S 2. Or we can obtain some compactification 
 using the embedding of $\pi$ to the projective space over $\Delta$ 
 in an obvious way. But we do not discuss
 about the singular fiber here because we do not use any geometric
 information of the singular fiber in the sequel. 
\end{remark}

Let $p:\mathcal A\to \Delta ^*$ be the family of abelian surfaces in 
Definition \ref{kumfam}. Let $\overline T_m$ 
be the associated monodromy operator on  
$H^m(\mathcal A_t)$ and assume that $\overline T_1$ is unipotent. Then
one can easily see that all $\overline T_m$'s are unipotent. 
We denote $\overline N_m=\log \overline T_m$. 
Of course $\nilp (\overline N_1)\leqslant 1$ and equality holds if
$\overline T_1$ is non-trivial. Assume $\nilp (\overline N_1)=1$ and let
$l=\rank \overline N_1\,(=1,2)$. Then we know that 
\begin{equation}\label{nilponas}
 \nilp (\overline N_2)=l, \quad \nilp (\overline N_3)=1, 
\quad \nilp (\overline N_4)=0.
\end{equation}
 
Take the family of generalized Kummer $2n$-folds 
$\pi :\Kum ^n(\mathcal A/\Delta ^*)\to \Delta ^*$ and 
let $T_m$ be the monodromy operator on $H^m(\Kum ^n(\mathcal A_t))$ and
$N_m=\log T_m$. As in \S 2, $N_2$ can be expressed by 
\[
 N_2=\overline N_2\oplus \id
\]
under $H^2(\Kum ^n(\mathcal A_t),\mathbb C)\cong H^2(A,\mathbb C)\oplus
\mathbb C\cdot [E_t]$. 
We calculate $\nilp (N_{2k})$ for $1<k\leqslant n$ below. 

We prepare a lemma.

\begin{lemma}\label{lowerbound}
 Let $\pi :\mathcal X\to \Delta$ be a degeneration of irreducible
 symplectic manifolds, and $T_m$ the associated monodromy operator 
 on $H^m(\mathcal X_t,\mathbb C)$.
 Assume $T_2$ and $T_{2k}$ are unipotent and take $N_m=\log T_m\;
 (m=2,2k)$. Then, $\nilp (N_{2k})\geqslant k\cdot \nilp (N_2)$. 
 In particular $N_{2k}^{lk}=0\; \Rightarrow\; N_2^l=0$.
\end{lemma}

\begin{proof}
 We again use the result of Verbitsky (Proposition \ref{Verb}). 
 We have an injective homomorphism $S^k=\Sym ^kH^2(\mathcal X_t)
 \hookrightarrow H^{2k}(\mathcal X_t)$ and we have
 \[
  \nilp (N_{2k})\geqslant \nilp (N_{2k|S^k})=\nilp (\log \Sym
 ^kT_2)=k\cdot\nilp (N_2)
 \]
 by Lemma \ref{tensornilp}.
\end{proof}

The main tool of our calculation is again 
the theorem of G\"ottsche--Soergel:

\begin{theorem}[G\"ottsche--Soergel \cite{GS}, 
 Theorem 7]\label{GSformula} 
 Let $A$ be an abelian surface and 
 fix a natural number $n$. Then, there is a canonical isomorphism
 \begin{equation}\label{GS}
 H^{i+2n}(A\times \Kum ^{n-1}(A))=\bigoplus _{\alpha\in P(n)} 
 \left( H^{i+2|\alpha |} (A^{(\alpha)}) \right)^{\oplus g(\alpha)^4},  
 \end{equation}
 where the notations $P(n)$, $|\alpha|$, and $A^{(\alpha)}$ are the same as
 in Theorem \ref{GSformula0}, and $g(\alpha )=\gcd \{k\mid \alpha _k\neq
 0\}$. 
\end{theorem}

Let us consider $\mathcal A\times _{\Delta ^*} \Kum ^n(\mathcal A/\Delta
^*)\to \Delta ^*$, which is in fact a degree $(n+1)^4$ 
\'etale cover of $\Hilb ^{n+1}(\mathcal A/\Delta ^*)$ over $\Delta ^*$, 
and the associated monodromy operator $\tilde T_m$ on 
$H^m(\mathcal A_t\times \Kum ^n(\mathcal A_t))$. 
Since G\"ottsche--Soergel decomposition \eqref{GS}
is induced by the action of
symmetric group $\mathfrak S_n$ on $\mathcal A^n\to \Delta ^*$
corresponding to $\Hilb ^n(\mathcal A/\Delta ^*)\to \Delta ^*$, 
the monodromy operator $\tilde T_m$ respects the decomposition
\eqref{GS}.

\begin{example}\label{kumn=2}
 Consider the case where $n=2$, i.e., the monodromies on the
 cohomologies  $H^m(\Kum ^2(\mathcal A_t))$ associated to $\pi :\Kum
 ^2(\mathcal A/\Delta ^*)\to \Delta ^*$.  
 We shall write $A=\mathcal A_t$ below for
 simplicity. By the G\"ottsche--Soergel formula (Theorem \ref{GSformula}), 
 we have
 \[
 H^4(A\times \Kum ^2(A))
 =H^0(A)^{\oplus 3^4}\oplus H^2(A\times A)\oplus H^4(A^{(3)})
 \]
 and one can easily check
 \[
 \begin{aligned}
  H^2(A\times A) &\cong H^2(A)^{\oplus 2}
  \oplus \left(H^1(A)\otimes H^1(A)\right), \\
  H^4(A^{(3)}) &\cong H^4(A)\oplus \left(H^3(A)\otimes H^1(A)\right) 
  \oplus \Sym ^2 H^2(A)\oplus \left(H^2(A)\otimes \wedge^2H^1(A)\right).
 \end{aligned}
 \]
 In particular, $\tilde T_m$ is unipotent. Take $\tilde N_m=\log \tilde
 T_m$. Using Lemma \ref{tensornilp} and \eqref{nilponas}, 
 we know that $\nilp (\tilde N_4)=2l$,
 where $l=\nilp (\overline N_2) \geqslant 1$. 
 On the other hand, the K\"unneth formula 
 \begin{multline*}
 H^4(A\times \Kum ^2(A))
  \cong H^4(\Kum ^2(A))
  \oplus \left(H^3(\Kum ^2(A))\otimes H^1(A)\right)\\
  \oplus \left(H^2(\Kum ^2(A))\otimes H^2(A)\right)
  \oplus H^4(A)
 \end{multline*}
 infers that $T_m$ is also unipotent and 
 \[
 \nilp (\tilde N_4)=\max\{
 \nilp (N_4), \nilp (N_3)+1, 2l\}
 \]
 if $l>0\; \Leftrightarrow \nilp (\overline N_1)=1$. 
 Therefore 
 \[
 2l=\nilp (\tilde N_4)\geqslant \nilp (N_4)\geqslant 2l, 
 \]
 where the last inequality is due to Lemma \ref{lowerbound}. 
 Note that this holds also for $l=\rank \overline N_1=0$. 
 In summary, we have 
 \begin{equation}\label{kumn=2res}
 \nilp (N_4)=2\cdot \nilp (N_2). 
 \end{equation}
\end{example}

Generalizing this method of calculation, we can show the following.

\begin{theorem}\label{kumn>2}
 Take $n\geqslant 2$ and $2\leqslant k\leqslant n$. 
 Let $T_m$ be the monodromy operator on 
 $H^m(\Kum ^n(\mathcal A_t),\mathbb C)$ associated with 
 the family $\pi :\Kum ^n(\mathcal A)\to \Delta ^*$ in Definition 
 \ref{kumfam}. Then, $T_m$ is also unipotent and 
 \begin{equation}\label{kumn>2res}
  \nilp (N_{2k})=kl,\quad \nilp (N_{2k-1})\leqslant kl-1.
 \end{equation}
 for $N_m=\log T_m$. In particular $\nilp (N_{2k})\in \{0,k,2k\}$. 
\end{theorem}

We prove this theorem via several steps. We keep the notation
$A=\mathcal A_t$. 

\begin{lemma}\label{kummer1}
 Let $a$ be a positive integer and 
 $N(m,a)$ the logarithm of the induced monodromy operator 
 on $H^m(A^{(a)})$. Then, 
 \[
  \begin{aligned}
   \nilp (N(2M,a)) &=\begin{cases}
		      Ml & (M\leqslant a)\\
		      (2a-M)l & (a<M\leqslant 2a)\\
		      0 & (M>2a)
		     \end{cases} ,\\
   \nilp (N(2M+1,a)) &=\begin{cases}
		      Ml+1 & (M\leqslant a)\\
		      (2a-M-1)l+1 & (a<M\leqslant 2a)\\
		      0 & (M>2a)
		     \end{cases} .
  \end{aligned}
 \]
\end{lemma}

\begin{proof}
 By the Poincar\'e duality, we have only to consider the case $M\leqslant
 a$. Let $\mu=(\mu_1,\cdots ,\mu_4)$ be the partition of $m=2M$ or
 $2M+1$
 \[
 m=\mu_1\cdot 1+\mu _2\cdot 2+\mu _3\cdot 3+\mu _4\cdot 4
 \]
 under $|\mu |=\mu _1+\mu _2+\mu _3+\mu _4\leqslant a$. Then we have 
 \[
 H^m(A^{(a)}) \cong \bigoplus _{\mu} H^{(\mu)}(A^{(a)})
 \]
 where
 \[ 
 H^{(\mu)}(A^{(a)})
 \cong
 \wedge ^{\mu _1}H^1(A)
 \otimes \Sym ^{\mu _2}H^2(A)\otimes \wedge ^{\mu _3} H^3(A)
 \otimes \Sym ^{\mu _4}H^4(A). 
 \]
 Let $N((\mu), a)$ be the logarithm of the induced monodromy operator 
 on $H^{(\mu)}(A^{(a)})$, then 
 $\displaystyle \nilp (N(m,a))=\max _{\mu} \nilp
 (N((\mu),a))$.  One can easily see that 
 $\nilp N((\mu), a)$ attains the maximum at 
 \[
 \begin{aligned}
  \mu &=(0,M,0,0) & (m=2M),\\
  \mu &=(1,M,0,0) & (m=2M+1), 
 \end{aligned}
 \]
 and the corresponding maximum values are $Ml$ and $Ml+1$, respectively. 
\end{proof}

\begin{lemma}\label{kummer2}
 Notation as above. $\nilp (\tilde N_{2k})=kl$.
\end{lemma}

\begin{proof}
 By the G\"ottsche--Soergel formula, 
 \[
 H^{2k}(A\times \Kum ^n(A))=\bigoplus _{\alpha \in P(n+1)}
 \left( H^{m_k(\alpha)}(A^{(\alpha)})\right)^{\oplus g(\alpha)^4}  
 \]
 where $m_k(\alpha)=2|\alpha |-2(n-k+1)$. 
 Let $\tilde N_{2k}(\alpha)$ be the logarithm of the induced monodromy 
 on $H^{m_k(\alpha)}(A^{(\alpha)})$ and $\nu(\alpha)=\# \{i\mid 
 \alpha _i\neq 0\}-1$. 
 Then, by the K\"unneth formula and
 Lemma \ref{kummer1}, we have
 \begin{multline*}
  \nilp (\tilde N_{2k}(\alpha)) 
  =\max\left\{\frac{m_k(\alpha)}2 \cdot l, 
  \left(\frac{m_k(\alpha)}2 -1\right)\cdot l+2,\right. \\ 
  \left. \cdots , \left(\frac{m_k(\alpha)}2 -\nu(\alpha)\right)\cdot l
  +2\nu (\alpha)\right\}.
 \end{multline*}
 As we have an obvious inequality $|\alpha |+\nu(\alpha)\leqslant n+1$, 
 which is immediate from $\alpha \in P(n+1)$, we have
 \[
  \left(\frac{m_k(\alpha)}2 -\nu (\alpha)\right)\cdot l+2\nu (\alpha)
 =\left(|\alpha|-n+k-\nu(\alpha)-1\right)\cdot l+2\nu (\alpha)
 \leqslant kl
 \]
 even if $l=1$. Therefore, we get 
 \[
  \nilp (\tilde N_{2k})=\max \{\nilp (\tilde N_{2k}(\alpha))
 \mid \alpha \in P(n+1)\}\leqslant
 kl. 
 \]
 But $\nilp (\tilde N_{2k}(\alpha))$ attains the maximum 
 $kl$ at $\alpha =(n+1,0,\cdots)$. This completes the proof of the
 lemma. 
\end{proof}

\begin{proof}[Proof of Theorem \ref{kumn>2}]
 By the K\"unneth formula
 \[
 H^{2k}(A\times \Kum ^n(A))
 \cong \bigoplus _{i=0}^4 H^i(A)\otimes H^{2k-i}(\Kum ^n(A)),
 \]
 we have
 \begin{multline}\label{kunnethestim}
 \nilp (\tilde N_{2k})
 =\max \{ \nilp (N_{2k}), \nilp (N_{2k-1})+1, \\ 
  \nilp (N_{2k-2})+l, \nilp (N_{2k-3})+1, \nilp (N_{2k-4})\}.
 \end{multline}
 In particular, $T_m$ is unipotent and 
 \[
 \nilp (N_{2k})\leqslant kl\,\quad \nilp (N_{2k-1})\leqslant kl-1
 \]
 by Lemma \ref{kummer2}. On the other hand, we have $\nilp
 (N_{2k})\geqslant kl$ by Lemma \ref{lowerbound}, so this proves the
 theorem.
\end{proof}

\begin{remark}
 We should note that the proof of 
 the equality \eqref{kumn>2res} 
 is limited to the case of the family of
 generalized Kummer varieties arising from a degeneration of abelian
 surfaces, unlike the case of Hilbert type (cf. Theorem \ref{MofHilbFam}).
\end{remark}

\section{Good degeneration of a compact symplectic K\"ahler manifold}

From what we have seen in the previous sections, we can easily pose
the following question: 
\emph{Do the monodromy operators on the cohomologies associated with the
degeneration of an irreducible symplectic manifold have some special
property?} To study the question in a somewhat general situation, it
is certainly one way to consider the relation between the monodromies and 
the geometry of the singular fiber. Along this direction, 
we have the powerful theory of the
limit mixed Hodge structure in the case of semi-stable degenerations.

\begin{definition}
 Let $M$ be a complex manifold. A divisor $D$ on $M$ is a \emph{simple
 normal crossing} divisor (SNC in short) if $D=\sum D_i$ is reduced, 
 every irreducible component $D_i$ is smooth, and for any point $p\in
 D$, the local equation of $D$ in $M$ is given by $x_0\cdots x_r=0$ for
 some $r$. A degeneration $\pi :\mathcal X\to \Delta$ is 
 \emph{semi-stable} if the total space $\mathcal X$ is smooth and 
 $X=\mathcal X_0=\pi ^{-1}(0)$ is a SNC divisor 
 as scheme theoretic fiber.
\end{definition}

By the semi-stable reduction theorem, a degenerating family is always
birational to a semi-stable one after taking some cyclic base change. 
It is also known that the monodromies $T_m$ on 
$H^m(\mathcal X_t,\mathbb C)$ is unipotent 
if the degeneration $\pi:\mathcal X\to \Delta$ 
is semi-stable. In this sense, we can consider a semi-stable degeneration
as a geometric counterpart of the concept of ``unipotent monodromy''.

But one should note that there are many semi-stable models for a given
degenerating family since one can operate birational modifications
keeping the family semi-stable. To carry out some geometric arguments on
the singular fiber, it is desirable to have a kind of ``minimality'' of
the family. As that kind of thing, we propose the following definition
of a good degeneration of a compact symplectic K\"ahler manifold.

\begin{definition}\label{gooddegsympl}
A \emph{good degeneration of compact symplectic K\"ahler manifold} is a 
degeneration $\pi : \mathcal{X}\to \Delta$ of 
relative dimension $2n$ satisfying
\begin{enumerate}[(i)]
\item $\pi$ is semi-stable.
\item There exists a relative logarithmic 2-form $\sigma _{\pi}\in 
      H^0(\mathcal X,\Omega ^2_{\mathcal X/\Delta}(\log X))$ such that 
      $\wedge ^n \sigma _{\pi} \in H^0(\mathcal X, K_{\mathcal X/\Delta})$ 
      is nowhere vanishing (see, for example, \cite{St} 
      for the definition of the logarithmic differential forms).
\end{enumerate}
\end{definition}

Note that the condition (ii) implies that $K_{\mathcal X/\Delta}$ is
trivial. In particular the definition above 
agrees with the definition of good degeneration of K3 surface, 
so-called Kulikov model, 
by Kulikov--Persson--Pinkham \cite{Ku,P-P} if $n=1$. 

Let us construct an example of a good degeneration of an irreducible
symplectic manifold using the degenerating family of the Hilbert schemes
on K3 surfaces in \S 2.

\begin{theorem}\label{exampleofdsv}
 Let $p : \mathcal S\to \Delta $ be a projective type II degeneration of K3 
 surface, i.e., $p$ is a projective good degeneration of K3 surface 
 with the singular fiber $\mathcal S_0 = S_0\cup S_1\cup \cdots 
 \cup S_{k-1} \cup S_k$ where $S_0$ and $S_k$ are rational surfaces, 
 $S_i\; (0< i < k)$ are elliptic ruled surfaces 
 and $S_i$ meets only $S_{i\pm 1}$ in smooth elliptic curves
 $C_i=S_i\cap S_{i+1}\; (i=0,\cdots k-1)$. 
 Consider the Hilbert scheme $\rho : \mathcal Y
 =\Hilb ^2(\mathcal S/\Delta)\to \Delta$ of relative sub-schemes of 
 length 2. Then there exists a projective birational morphism 
 $\mu :\mathcal X \to \mathcal Y$ such that 
 \[
 \pi =\rho \circ \mu : \mathcal X \to \Delta
 \] 
 is a good degeneration of compact symplectic K\"ahler manifold.
\end{theorem}

In the proof of the theorem given below, it is essential 
the condition that the length of sub-schemes in question is 2. 
It is natural to ask either the same conclusion holds for the larger 
length sub-schemes. To answer this question, it is likely that 
a more intrinsic interpretation of the resolved space $\mathcal X$ is needed.

\begin{proof}
 Note that we can assume that $\mathcal Y$ is a K\"ahler
 space (Example \ref{HilbFam}) so that $\mathcal X$ in the theorem is
 automatically a K\"ahler manifold as far as $\mu$ is a projective
 resolution.

The family $p : \mathcal S\to \Delta$ induces the morphism
\[
f: \Hilb ^2(\mathcal S)\to \Sym ^2(\mathcal S)\to \Sym ^2 (\Delta) .
\]
We can consider 
$\Hilb ^2(\mathcal S/\Delta)$ as a closed sub-scheme of
$\Hilb ^2(\mathcal S)$. Moreover if we define 
$d : \Delta \to \Sym ^2 (\Delta )$ by \linebreak[2]
$z\mapsto 2[z]$, then 
we have the commutative diagram 
\[
\xymatrix{
\Hilb ^2(\mathcal S / \Delta )\ar[r]\ar[dr]_{\rho}
 & \Hilb ^2(\mathcal S) \times
 	_{\Sym^2(\Delta )}\Delta \ar[d] \ar[r]
 & \Hilb ^2(\mathcal S)\ar[d]^f\\
 & \Delta \ar[r]^d & \Sym ^2(\Delta)
}\lower 47pt \hbox{.}
\]
A singular sub-scheme 
of length 2 is given as the base point plus tangent 
direction at that point and any sub-scheme of length 2 is 
given as a limit of smooth sub-schemes (i.e., sub-schemes 
consisting of different 2 points). 
Therefore 
\[
\Hilb ^2(\mathcal S)=
Bl _{D}(\mathcal S\times \mathcal S)/\mathfrak S_2
\]
where $D$ is the diagonal, $Bl _{D}$ stands for the blowing up 
along $D$ and $\mathfrak S_2$ acts as the permutation of components. 
Moreover $\Hilb ^2(\mathcal S/\Delta )$ is the closure of the open subset 
of $\Hilb ^2(\mathcal S) \times _{\Sym^2(\Delta )}\Delta$ 
consisting of points corresponding to the smooth sub-schemes. 
Thus we have the following commutative diagram
\[
\xymatrix{
 Bl _{D}(\mathcal S\times \mathcal S)\ar[r]^g\ar[d]
 & \mathcal S\times \mathcal S \ar[r]^h\ar[d]
 & \Delta \times \Delta \ar[d]\\
\Hilb ^2(\mathcal S)\ar[r]
 & \Sym ^2(\mathcal S)\ar[r]
 & \Sym ^2(\Delta)\\
\Hilb ^2(\mathcal S/\Delta)\ar[u]\ar[rr]^{\rho}
 && \Delta \ar[u]
}\lower 84pt\hbox{.}
\]
Let $W$ be the strict transform 
on $Bl _{D}(\mathcal S\times \mathcal S)$
of the inverse image by $h$ of
the diagonal of $\Delta\times \Delta$.
Then we have the commutative diagram 
\[
\xymatrix{
W\ar[r]^{(h\circ g)_{|W}}\ar[d] & \Delta\ar@{=}[d]\\
\mathcal Y=\Hilb(\mathcal S/\Delta )\ar[r]^{\quad\qquad\rho} & \Delta
}
\]
where the first vertical arrow is the quotient map by the action of 
$\mathfrak S_2$. 

 Now, since $\mathcal S$ is a type II degeneration of K3 surfaces, 
 components of $\mathcal Y_0=\rho ^{-1}(0)$ consists of 
 $Y_{ij}\; (0\leqslant i\leqslant j\leqslant k)$ where 
 $Y_{ii}=\Hilb ^2(S_i)$ and $Y_{ij}$ is isomorphic to 
 $S_i\times S_j$ if $j>i+1$. $Y_{i,i+1}$ is a singular variety given by 
 identifying the points in $C_i\times C_i\subset S_i\times S_{i+1}$ by
 the natural action of $\mathfrak S_2$ on $C_i\times C_i$. 
 The configuration of these components is as follows.
 \begin{center}
\unitlength 0.1in
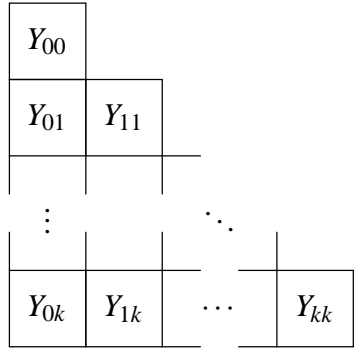
\begin{picture}( 19.6000, 18.0000)(  2.4000,-20.0000)
%
\special{pn 8}%
\special{pa 400 200}%
\special{pa 800 200}%
\special{pa 800 600}%
\special{pa 400 600}%
\special{pa 400 200}%
\special{fp}%
%
\special{pn 8}%
\special{pa 400 600}%
\special{pa 800 600}%
\special{pa 800 1000}%
\special{pa 400 1000}%
\special{pa 400 600}%
\special{fp}%
%
\special{pn 8}%
\special{pa 800 600}%
\special{pa 1200 600}%
\special{pa 1200 1000}%
\special{pa 800 1000}%
\special{pa 800 600}%
\special{fp}%
%
\special{pn 8}%
\special{pa 400 1600}%
\special{pa 800 1600}%
\special{pa 800 2000}%
\special{pa 400 2000}%
\special{pa 400 1600}%
\special{fp}%
%
\special{pn 8}%
\special{pa 800 2000}%
\special{pa 1200 2000}%
\special{pa 1200 1600}%
\special{pa 800 1600}%
\special{pa 800 2000}%
\special{fp}%
%
\special{pn 8}%
\special{pa 1800 1600}%
\special{pa 2200 1600}%
\special{pa 2200 2000}%
\special{pa 1800 2000}%
\special{pa 1800 1600}%
\special{fp}%
\put(6.0000,-4.0000){\makebox(0,0){$Y_{00}$}}%
\put(6.0000,-8.0000){\makebox(0,0){$Y_{01}$}}%
\put(10.0000,-8.0000){\makebox(0,0){$Y_{11}$}}%
\put(6.0000,-13.0000){\makebox(0,0){$\vdots$}}%
\put(15.0000,-13.0000){\makebox(0,0){$\ddots$}}%
\put(15.0000,-18.0000){\makebox(0,0){$\cdots$}}%
\put(6.0000,-18.0000){\makebox(0,0){$Y_{0k}$}}%
\put(10.0000,-18.0000){\makebox(0,0){$Y_{1k}$}}%
\put(20.0000,-18.0000){\makebox(0,0){$Y_{kk}$}}%
%
\special{pn 8}%
\special{pa 400 1000}%
\special{pa 400 1200}%
\special{fp}%
\special{pa 800 1200}%
\special{pa 800 1000}%
\special{fp}%
\special{pa 1200 1000}%
\special{pa 1200 1200}%
\special{fp}%
\special{pa 1400 1000}%
\special{pa 1200 1000}%
\special{fp}%
\special{pa 400 1600}%
\special{pa 400 1400}%
\special{fp}%
\special{pa 800 1600}%
\special{pa 800 1400}%
\special{fp}%
\special{pa 1200 1400}%
\special{pa 1200 1600}%
\special{fp}%
\special{pa 1200 1600}%
\special{pa 1400 1600}%
\special{fp}%
\special{pa 1400 2000}%
\special{pa 1200 2000}%
\special{fp}%
\special{pa 1800 2000}%
\special{pa 1600 2000}%
\special{fp}%
\special{pa 1600 1600}%
\special{pa 1800 1600}%
\special{fp}%
\special{pa 1800 1600}%
\special{pa 1800 1400}%
\special{fp}%
\end{picture}%
 \end{center}
 Now we define $\mu :\mathcal X\to \mathcal Y$ as the blowing up 
 along $\displaystyle 
 C = \hskip -8pt
 \bigcup _{0\leqslant 2s \leqslant 2t \leqslant k} 
 \hskip -8pt Y_{2s,2k}$ and show that $\mu$ is a small resolution 
 and $\pi=\rho\circ\mu :\mathcal X\to \Delta$ is semi-stable.

 Noting that $\rho$ is smooth over 
 $\Delta ^*=\Delta -\{0\}$, we have only to look at the singular fiber.
 Let $Z$ be a sub-scheme of length 2 on $\mathcal S_0$. 
 If $\supp (Z)$ is contained in the smooth locus of $\mathcal S_0$, 
 $\rho$ is obviously smooth at the point $[Z]\in \mathcal Y$. 
 Suppose that exactly two components, say $Y_{ij}$ and $Y_{i,j+1}$, 
 meet at $[Z]$.  
 Then $Z=q_1+q_2$ where $q_1\in S_i\backslash \cup _{k\neq i} S_k$ 
 and $q_2\in S_j \cap S_{j+1}$. Locally at $[Z]$, 
 $\Hilb ^2(\mathcal S)$ is a product of open neighborhoods of $q_1$ and
 $q_2$ since $[Z]$ is away from diagonal. Hence, $f$ is given by 
 \[
  (z_1,\dots ,z_6)\mapsto (z_1,z_4z_5)
 \]
 for some appropriate coordinate. $\mathcal Y$ is defined by 
 $z_1-z_4z_5=0$ and $\mathcal Y_0$ is given by $z_1=z_4z_5=0$. 
 Therefore $\mathcal Y$ have a coordinate $(z_2, \dots, z_6)$ at $[Z]$
 such that $\rho$ is given by $(z_2, \dots, z_6)\mapsto z_4z_5$. This
 shows that $\mathcal Y$ is non-singular and $\rho$ is semi-stable at
 $[Z]$. Therefore, if $\rho$ is not semi-stable at $[Z]$, the $\supp
 (Z)$ must be contained in $\Sing (\mathcal S_0)$. In other words, the
 points where $\rho$ is not semi-stable are 
 contained in the locus of points on $\mathcal Y_0$ where at least 3
 components of $\mathcal Y_0$ intersect. 

Now assume $\supp (Z)\subset \Sing (\mathcal S_0)$. 
Away from the diagonal, $\Hilb ^2(\mathcal S)$ is 
locally the direct product of open neighborhoods of $\mathcal S$ 
and $f$ is given by
\[
(z_1,\dots ,z_6)\mapsto (z_1z_2,z_4z_5)
\]
with respect to some coordinate at $[Z]$ 
so that total space of $\rho$ is defined by $z_1z_2-z_4z_5=0$ 
and the central fiber of $\rho$ is defined by $z_1z_2=z_4z_5=0$. 
As the center $C$ of $\mu$ is defined by, say, $z_1=z_4=0$, 
the defining equation of $\mathcal X$ is given by
$z'_1z'_2-z'_5=0$ and $\pi=\rho\circ \mu$ is described by
\[
(z'_1,\dots ,z'_6)\mapsto z'_1z'_2=z'_5
\]
This shows $\mathcal X$ is smooth and $X=\pi ^{-1}(0)$ is normal 
crossing divisor in this coordinate neighborhood. 

Next we consider local description of $[Z]\in \rho ^{-1}(0)$ with 
$\supp Z=\{q\},\; q\in \Sing (\mathcal S_0)$. Take a local coordinate 
$(z_1,\dots ,z_6)$ of $(q,q)\in \mathcal S\times \mathcal S$. 
The defining equation of diagonal $D$ is $z_1-z_4=z_2-z_5=z_3-z_6=0$. 
Put 
\begin{eqnarray*}
x_1=z_1+z_4, & x_3=z_2+z_5, & x_5=z_3+z_6,\\
x_2=z_1-z_4, & x_4=z_2-z_5, & x_6=z_3-z_6.
\end{eqnarray*}
The map $h$ is given by
\[
(x_1,\dots ,x_6)\mapsto ((x_1+x_2)\cdot (x_3+x_4),
(x_1-x_2)\cdot (x_3-x_4)).
\]
The defining ideal of $D$ is $(x_2,x_4,x_6)$ and a piece of 
the blowing up $g$ is locally described by
\begin{eqnarray*}
x_1=y_1, & x_2=y_2y_6, & x_3=y_3, \\
x_4=y_4y_6, & x_5=y_5, & x_6=y_6,
\end{eqnarray*}
and $h\circ g$ is given by 
\[
(y_1,\dots ,y_6)\mapsto 
((y_1+y_2y_6)(y_3+y_4y_6),(y_1-y_2y_6)(y_3-y_4y_6)).
\]
Restricting on $W$, i.e., imposing the condition
\[
y_2y_3+y_1y_4=0,
\]
$h\circ g$ is given by
\[
(y_1,\dots ,y_6)\mapsto 
y_1y_3+y_2y_4y_6^2 .
\]
Noting that the action of $\mathfrak S_2$ is
 $\mbox{diag}(1,1,1,1,1,-1)$, $\rho$ is described by
\[
(w_1,\dots ,w_6)\mapsto w_1w_3+w_2w_4w_6
\]
under the condition $w_2w_3+w_1w_4=0$. In this coordinate, 
the equation of the center $C$ of $\mu$ is, say, $w_1=w_2=0$. 
Therefore the equation of $\mathcal X$ is given by 
$w'_3+w'_1w'_4=0$ and $\pi =\rho\circ\mu$ is described by
\[
(w'_1,\dots ,w'_6)\mapsto w'_1w'_2w'_3+w'_2w'_4w'_6 .
\]
The fiber $\pi ^{-1}(t)$ is $w'_2w'_4(w'_6-w_1^{\prime\; 2})=t$
in the coordinate $(w'_1, w'_2, w'_4, w'_5, w'_6)$. 
These calculations show that $\pi $ is semi-stable also in 
this coordinate neighborhood. 
 
 In summery, we proved that 
 $\mu$ is a small resolution and $\rho$ is everywhere semi-stable. 
 We remark that the configuration of
 components of $\mathcal X_0=\pi ^{-1}(0)$ is as follows.
 \begin{center}
  \unitlength 0.1in
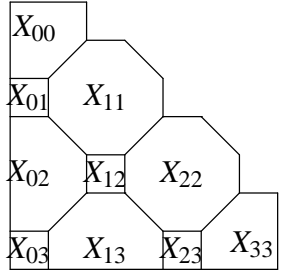
  \begin{picture}( 16.6000, 14.0000)(  1.4000,-16.0000)
   %
   \special{pn 8}%
   \special{pa 400 200}%
   \special{pa 400 600}%
   \special{fp}%
   \special{pa 400 600}%
   \special{pa 600 600}%
   \special{fp}%
   \special{pa 600 600}%
   \special{pa 800 400}%
   \special{fp}%
   \special{pa 800 400}%
   \special{pa 800 200}%
   \special{fp}%
   \special{pa 800 200}%
   \special{pa 400 200}%
   \special{fp}%
   %
   \special{pn 8}%
   \special{pa 600 600}%
   \special{pa 600 800}%
   \special{fp}%
   \special{pa 600 800}%
   \special{pa 800 1000}%
   \special{fp}%
   \special{pa 800 1000}%
   \special{pa 1000 1000}%
   \special{fp}%
   \special{pa 1000 1000}%
   \special{pa 1200 800}%
   \special{fp}%
   \special{pa 1200 800}%
   \special{pa 1200 600}%
   \special{fp}%
   \special{pa 1200 600}%
   \special{pa 1000 400}%
   \special{fp}%
   \special{pa 1000 400}%
   \special{pa 800 400}%
   \special{fp}%
   %
   \special{pn 8}%
   \special{pa 400 600}%
   \special{pa 400 800}%
   \special{fp}%
   \special{pa 400 800}%
   \special{pa 600 800}%
   \special{fp}%
   %
   \special{pn 8}%
   \special{pa 400 800}%
   \special{pa 400 1400}%
   \special{fp}%
   \special{pa 400 1400}%
   \special{pa 600 1400}%
   \special{fp}%
   \special{pa 600 1400}%
   \special{pa 800 1200}%
   \special{fp}%
   \special{pa 800 1200}%
   \special{pa 800 1000}%
   \special{fp}%
   %
   \special{pn 8}%
   \special{pa 800 1200}%
   \special{pa 800 1200}%
   \special{fp}%
   \special{pa 1000 1200}%
   \special{pa 1000 1200}%
   \special{fp}%
   %
   \special{pn 8}%
   \special{pa 800 1200}%
   \special{pa 1000 1200}%
   \special{fp}%
   \special{pa 1000 1200}%
   \special{pa 1000 1000}%
   \special{fp}%
   %
   \special{pn 8}%
   \special{pa 400 1400}%
   \special{pa 400 1600}%
   \special{fp}%
   \special{pa 400 1600}%
   \special{pa 600 1600}%
   \special{fp}%
   \special{pa 600 1600}%
   \special{pa 600 1400}%
   \special{fp}%
   %
   \special{pn 8}%
   \special{pa 1000 1200}%
   \special{pa 1200 1400}%
   \special{fp}%
   \special{pa 1200 1400}%
   \special{pa 1200 1600}%
   \special{fp}%
   \special{pa 1200 1600}%
   \special{pa 600 1600}%
   \special{fp}%
   %
   \special{pn 8}%
   \special{pa 1200 800}%
   \special{pa 1200 800}%
   \special{fp}%
   \special{pa 1400 800}%
   \special{pa 1400 800}%
   \special{fp}%
   \special{pa 1200 800}%
   \special{pa 1400 800}%
   \special{fp}%
   \special{pa 1400 800}%
   \special{pa 1600 1000}%
   \special{fp}%
   \special{pa 1600 1000}%
   \special{pa 1600 1200}%
   \special{fp}%
   \special{pa 1600 1200}%
   \special{pa 1400 1400}%
   \special{fp}%
   \special{pa 1400 1400}%
   \special{pa 1200 1400}%
   \special{fp}%
   %
   \special{pn 8}%
   \special{pa 1200 1600}%
   \special{pa 1400 1600}%
   \special{fp}%
   \special{pa 1400 1600}%
   \special{pa 1400 1400}%
   \special{fp}%
   %
   \special{pn 8}%
   \special{pa 1400 1600}%
   \special{pa 1800 1600}%
   \special{fp}%
   \special{pa 1800 1600}%
   \special{pa 1800 1200}%
   \special{fp}%
   \special{pa 1800 1200}%
   \special{pa 1600 1200}%
   \special{fp}%
   \put(9.0000,-7.0000){\makebox(0,0){$X_{11}$}}%
   \put(9.0000,-11.0000){\makebox(0,0){$X_{12}$}}%
   \put(13.0000,-11.0000){\makebox(0,0){$X_{22}$}}%
   \put(5.0000,-15.0000){\makebox(0,0){$X_{03}$}}%
   \put(9.0000,-15.0000){\makebox(0,0){$X_{13}$}}%
   \put(13.0000,-15.0000){\makebox(0,0){$X_{23}$}}%
   \put(16.7000,-14.7000){\makebox(0,0){$X_{33}$}}%
   \put(5.0000,-7.0000){\makebox(0,0){$X_{01}$}}%
   \put(5.0000,-11.0000){\makebox(0,0){$X_{02}$}}%
   \put(5.3000,-3.4000){\makebox(0,0){$X_{00}$}}%
  \end{picture}%
  
  This is the case where $k=3$.
 \end{center}
 
 It remains to show that there exists a relative logarithmic  
 symplectic form $\sigma _{\pi} \in 
 H^0(\mathcal X,\Omega ^2_{\mathcal X/\Delta}(\log X))$.
 The original family $p : \mathcal S \to \Delta$ have nowhere vanishing
 relative logarithmic 2-form $\omega$. We claim that the relative 2-form 
 $\widetilde{\omega} =\mbox{pr}_1 ^* \omega +\mbox{pr}_2 ^*\omega $ on
 $\mathcal S\times \mathcal S$ induces $\sigma _{\pi}$. 
 The restriction of $\widetilde\omega$ to 
 the inverse image by $h$ of the diagonal of $\Delta \times \Delta$ 
 is invariant under the action of $\mathfrak S_2$, therefore 
 $(g^*\widetilde\omega )_{|W}$ descends to 
 a relative log 2-form $\varphi$ on 
 $\Hilb ^2(\mathcal S/\Delta )$. Let $\sigma _{\pi} = \mu ^* 
 \varphi$. As $\varphi$ is non-degenerate outside the critical
 points of $\pi$ by the argument same as in \cite{Be}, we know that
 $K_{\mathcal Y}=0$.  By the theory of limiting Hodge structure 
 induced by a semi-stable degeneration \cite{St}, the direct image sheaf
 $\pi _*\Omega ^2_{\mathcal X/\Delta}(\log X)$ is locally free. 
 Since $\mu :\mathcal X\to \mathcal Y$ is a small resolution, 
 $K_{\mathcal X}$ is trivial. Therefore 
 we see that $\sigma _{\pi}$ above defines everywhere non-degenerate 
 section of the invertible sheaf 
 $\pi _*\Omega ^2_{\mathcal X/\Delta}(\log X)$. 
 This shows that $\sigma _{\pi}$ satisfies the condition 
 of Definition \ref{gooddegsympl}, (ii).
\end{proof}

\begin{remark}
 We may construct similar example from type III degenerations of 
 K3 surfaces. But in these cases, the combinatorics of the components 
 becomes more complicated and the choice of the center of a resolution
 should be subtle.
\end{remark}

We should also remark that our definition of good degeneration 
may be ``too good'' in general. 
In view of the situation of the complexity of 
the minimal models in higher dimensions, it is too optimistic to expect 
a good degeneration model for a given degeneration of 
irreducible symplectic K\"ahler manifolds. For example,
Kulikov--Persson--Pinkham model can be seen 
as a relatively minimal 3-fold model 
over the base and their construction of the
good model is already quite complicated.

\section{Monodromy of a good degeneration of an irreducible symplectic manifold}

Now we consider the behavior of the monodromy operators on the
cohomologies associated to a good degeneration of compact irreducible
symplectic K\"ahler manifolds. 

Inspired by Theorems \ref{MofHilbFam} and \ref{kumn>2}, 
we propose the following conjecture.

\begin{conjecture}\label{conj1}
 Let $\pi:\mathcal X\to \Delta$ be a degeneration of an irreducible
 symplectic $2n$-fold, $T_m$ the associated monodromy operator on 
 $H^m(\mathcal X_t,\mathbb C)$. Put 
 $N_m=\log T_m$ and assume $T_{2k}$ is unipotent for $k\leqslant n$. 
 Then $\nilp (N_{2k})=k\cdot \nilp (N_2)$ for $k\leqslant n$. 
\end{conjecture}

We can also consider a weak version of this conjecture.

\begin{conjecture}\label{conj2}
 Under the assumption and notation as in Conjecture \ref{conj1}, 
 we have $\nilp (N_{2k})\in\{0,k,2k\}$. 
\end{conjecture}

Of course, Conjecture \ref{conj1} implies Conjecture \ref{conj2}. 
We have already seen that these conjectures are true if 
\begin{enumerate}[(i)]
 \item $\pi$ is a degeneration of $\Hilb ^n(S)$, where $S$ is a K3 surface
       (Theorem \ref{MofHilbFam}). 
 \item $\pi$ is the family of generalized Kummer manifolds 
       $\pi :\Kum ^n(\mathcal A/\Delta)\to \Delta$ as in Definition
       \ref{kumfam} (Theorem \ref{kumn>2}). 
\end{enumerate}
These conjectures could be too naive for higher dimensions.
The author suspect that the conjectures may 
be true at least for lower dimensions, for example $2n=4$ and $6$.

Under the assumption of the existence of a good degeneration model, 
we can easily prove the following theorem, 
which is a partial answer to Conjecture \ref{conj2}. 

\begin{theorem}\label{main}
 Let $\pi:\mathcal X\to \Delta$ be a good degeneration of irreducible 
 symplectic $2n$-folds. Let $H_t^m=H^m(\mathcal X_t,\mathbb C)\; (t\neq
 0)$, $T_m$ be the monodromy operator on $H^m_t$ associated to the
 family $\pi$, and $N_m=\log T_m$. Take $k\leqslant n$ and 
 assume $N_{2k}^k=0$, then
 $N_{2k}=0$. In other words, $\nilp(N_{2k})\in\{0, k, k+1, \cdots, 2k\}$.
\end{theorem}

First we prepare a basic lemma.

\begin{lemma}\label{HD}
 Let $\pi: \mathcal X\to \Delta$ be a semi-stable degeneration, 
 $X=\mathcal X_0=\pi^{-1}(0)$ the singular fiber, 
 $H^m=H^m(X,\mathbb C)$ and $H^m_t=H^m(\mathcal X_t,\mathbb C)$ for
 $t\neq 0$. Consider the sheaf of torsion free differentials 
 $\hat \Omega ^p_X=\Omega ^p_X/(\mbox{Torsion})$ and the sheaf of
 logarithmic differentials 
 $\Omega ^m_X(\log )=\Omega ^m_{\mathcal X/\Delta}(\log X)\otimes
 \mathscr O_X$. 
 Then the $F^m$-part \[
  F^m\Psi : F^mH^m\to F^mH^m_t
 \]
 of the natural morphism $\Psi :H^m\to H^m_t$ of the 
 mixed Hodge structures is nothing but a map
 $H^0(X,\hat \Omega ^m_X)\to H^0(X,\Omega ^m_X(\log ))$ induced by the 
 natural inclusion $\hat \Omega ^m_X\to \Omega ^m_X(\log )$.
\end{lemma}

\begin{proof}
 First, we recall the fact that there is a resolution
 $\mathbb C_X\overset{\varepsilon}{\lto}\hat \Omega^{\bullet}_X$ 
 of the constant sheaf and the induced Hodge spectral sequence 
 \[
 {}_F E^{pq}_1=H^q(X,\hat \Omega ^p_X)\Longrightarrow H^{p+q}(X,\mathbb C)
 \]
 is $E_1$-degenerate (\cite{Fr1}, Proposition 1.5).    
 Moreover the associated filtration is the Hodge filtration 
 of the standard mixed Hodge structure on $H^m(X,\mathbb C)$. 

 Next, let us consider the real blow-up $\tilde\rho :
 \widetilde X\to X$ as in \cite{KN} \S 4. Let $d$ be the dimension of
 $X$ and $V_{r,d}=\{(z_0,\cdots, z_d)\in \mathbb C^{d+1} \mid z_0 \dots
 z_r=0\}$. Define the real blow-up $\rho:\widetilde V_{r,d} \to V_{r,d}$
 by
 \begin{multline*}
 \widetilde V_{r,d}=\{(s_0,\theta _0,\dots, s_r,\theta _r, z_{r+1},
 \dots ,z_d)
 \in (\mathbb R_{\geqslant 0}\times S^1)^{r+1}\times \mathbb C^{d-r}\\
 \mid s_0\cdots s_r=0,\; \theta _0+\dots +\theta _r=0 \}
 \end{multline*}
 and the obvious relation $z_j=s_j\cdot e^{\sqrt{-1}\, \theta _j}\; 
 (0\leqslant j\leqslant r)$. 
 One can easily check that this local construction is compatible with
 coordinate change so that we have the global real blow-up
 $\rho:\widetilde X\to X$ (\cite{KN}, p.404-405). 

 We can regard $\rho$ homotopically as the collapsing map 
 $\mathscr C_t:\mathcal X_t\to X$, which is the restriction of 
 the Clemens' retraction map $\mathscr C:\mathcal X\to X$ defined by
 Theorem 6.9 in \cite{C} (reader can find a more readable and short 
 summery in \cite{LTY}, \S 4). The local description of the collapsing
 map $\mathscr C_t$ is given as following; first we define a continuous
 map 
 \[
 \rE \mathscr C_t:\{s'=(s'_0,\dots ,s'_r)\in \mathbb R^{r+1} \mid 
 s'_0\cdots s'_r=t \}\to \{s_0\cdots s_r=0\}
 \]
 for small enough real number $t\neq 0$, 
 which is a homotopy equivalence (for the construction, see \cite{LTY}, 
 \S 4, or the original \cite{C}, \S 6). Then, define
 \[
 \mathscr C_t: \{(w_0,\cdots ,w_d)\in \mathbb C^{d+1}\mid w_0\cdots w_r=t\}
 \to V_{r,d}
 \]
 by
 \begin{multline*}
  (w_0,\cdots ,w_r; w_{r+1},\cdots ,w_d) \\
  \mapsto 
  ((\rE \mathscr C_t(s'))_0\cdot e^{\sqrt{-1}\theta _0},\cdots,
  (\rE \mathscr C_t(s'))_r\cdot e^{\sqrt{-1}\theta _r};
  w_{r+1},\cdots ,w_d), 
 \end{multline*}
 where $w_j=s'_j\cdot e^{\sqrt{-1}\theta _j}\; (0\leqslant j \leqslant
 r)$. 
 But this map factors
 as $\mathcal X_t\overset{\nu}{\lto}\widetilde X
 \overset{\tilde\rho}{\lto} X$ by
 \begin{multline*}
  (w_0,\cdots ,w_r; w_{r+1},\cdots ,w_d)\\
  \overset{\nu}{\mapsto} 
 ((\rE \mathscr C_t(s'))_0,\theta _0,\cdots,
 (\rE \mathscr C_t(s'))_r, \theta _r;
  w_{r+1},\cdots ,w_d). 
 \end{multline*}
 The map $\nu$ is obviously a homotopy equivalence 
 so that it induces an
 isomorphism $H^m(\widetilde X,\mathbb C)\overset{\sim}{\to}
 H^m_t$. Using this factorization, we get a decomposition of $\Psi$ as 
 \[
 H^m\overset{\tilde\rho^*}{\lto} H^m(\widetilde X,\mathbb C)
 \overset{\nu ^*}{\lto} H^m_t
 \]
 where the second arrow is an isomorphism. 

 Moreover, we have a quasi-isomorphism 
 $\tilde\varepsilon : \mathbb R\tilde\rho _*\mathbb C_{\widetilde X}
 \to \Omega ^{\bullet}_X(\log )$ (\cite{KN}, p.405) and 
 the induced Hodge spectral sequence
 \[
 {}_F E^{pq}_1=H^q(X,\Omega ^p_X(\log ))\Longrightarrow
 H^{p+q}(\widetilde X,\mathbb C)\cong H^{p+q}_t
 \]
 is $E_1$-degenerate (\cite{KN}, Lemma 4.1). 
 The associated filtration is the Hodge filtration of
 the limit Hodge structure on $H^m_t$ (See \cite{St}. A reader can find
 a readable survey by Zucker in \cite{Topics}, Chapter VII).
 From the construction of $\tilde\varepsilon$, we have the following 
 commutative diagram
 \[
 \xymatrix{
 \mathbb C_X\ar[d]^{\varepsilon} \ar[r]^{j\quad} 
  & \mathbb R\tilde\rho _*\mathbb C_{\tilde X} \ar[d]^{\tilde\varepsilon}\\
 \hat\Omega ^{\bullet}_X \ar[r]^{i\quad} & \Omega ^{\bullet}_X(\log )
 }
 \]
 where $i$ and $j$ are the natural morphisms. The lemma follows
 from the fact that the map
 $\Psi :H^m\to H^m_t$ is the same thing as the morphism induced by $j$ 
 and $\Psi$ is a morphism of mixed Hodge structures of weight $(0,0)$.
\end{proof}

\begin{remark}
 We can actually modify the definition of the collapsing map 
 $\mathscr C_t$ in a way that $\mathscr C_t$ is 
 a homeomorphism as is explained in \cite{P}.
\end{remark}

\begin{proof}[Proof of the Theorem]
 Let $H^m=H^m(X,\mathbb C)$ be the cohomology of the singular fiber
 $X=\mathcal X_0$ as in the previous lemma. 
 Consider the Clemens--Schmid exact sequence (a convenient reference is
 Morrison's lecture note, \cite{Topics}, Chapter VI):
 \begin{equation}\label{CS}
 F^m\Gr_m^WH^m\overset{\overline\Psi}{\lto}
 F^m\Gr_m^WH_t^m\overset{N_m}{\lto}F^{m-1}\Gr_{m-2}^WH^m_t. 
 \end{equation}
 The last term must be zero because $\Gr_{m-2}^WH^m_t$ is a pure Hodge
 structure of weight $(m-2)$. 
 
 As $\mathcal X_t$ is an irreducible symplectic manifold, 
 we have $h^{2,0}(\mathcal X_t)=1$. 
 The condition $N_{2k}^k=0$ implies 
 $N_2=0$ by Lemma \ref{lowerbound}, hence $\Gr _3^WH^2_t=\Gr
 _4^WH^2_t=0$. Therefore $F^2\Gr_2^WH^2_t$ is isomorphic to 
 $F^2H^2_t=H^0(X,\Omega ^2 _X(\log ))$. 
 By Lemma \ref{HD}, we have $F^2H^2=H^0(\hat \Omega ^2_X)$
 and $F^2\Psi$ is identified with the natural morphism 
 $H^0(X,\hat \Omega ^2_X)\to H^0(X,\Omega ^2_X(\log ))$. 
 Since we assumed that $\pi:\mathcal X\to \Delta$ is a good
 degeneration, we have a logarithmic symplectic form $\sigma \in
 H^0(X,\Omega ^2_X(\log ))$. Since $\overline\Psi$ is surjective, 
 this $\sigma$ lifts to an element of
 $H^0(X,\hat \Omega ^2_X)$. 

 Let $X^{[p]}=\bigsqcup X_{i_0\cdots i_p}$
 where $X_{i_0\cdots i_p}=X_{i_0}\cap \cdots \cap X_{i_p}$
 the $p$-fold intersection of the components of the singular fiber 
 $X=\sum X_i$. 
 Then, the weight spectral sequence
 \begin{equation} \label{SNCweight}
 {}_WE^{p\, q}_1=H^q(X^{[p]},\mathbb C)\Rightarrow
 E^{p+q}=H^{q+p}(X,\mathbb C)
 \end{equation}
 is $E_2$-degenerate. Assume $N_{2k}\neq 0$. 
 Then $X^{[1]}\neq \emptyset$ by \eqref{SNCweight} 
 and Clemens--Schmid exact sequence. 
 $\sigma\in H^0(X,\hat\Omega ^2_X)$ implies that $\sigma ^{\wedge n}$
 vanishes as a section of the canonical sheaf $\omega _X$ 
 at the generic points of the image of $X^{[1]}$ on $X$. This
 contradicts to non-degeneracy of the log symplectic form $\sigma$.
\end{proof}

\begin{remark}
 The proof shows in fact that $N_2=0 \Leftrightarrow N_{2k}=0$. Therefore,
 this theorem is also a partial answer to Conjecture \ref{conj1}. 
\end{remark}

\begin{corollary}\label{nocycle}
 Notation as above. For a good degeneration of irreducible
 symplectic $2n$-folds with non-trivial monodromy on the middle cohomology, 
 we have $X^{[p]}\neq \emptyset$ for $p\leqslant
 n$. In other words,
 for the dual graph $\Gamma$ of the configuration of the irreducible
 components of the singular fiber $X$, the dimension of the topological
 realization $|\Gamma |$ is at least $n$. 
\end{corollary}

\begin{proof}
 According to the weight spectral sequence \eqref{SNCweight}, we have 
 \begin{equation}\label{sncgr}
 \Gr_q^WH^{2n}
 \cong \frac{\Ker \big( H^q(X^{[2n-q]})\to H^q(X^{[2n-q+1]})\big)}
 {\iM \big( H^q(X^{[2n-q-1]})\to H^q(X^{[2n-q]})\big)}.
 \end{equation}
 By the Clemens--Schmid exact sequence, we also have 
 \[
 \Gr_q^WH^{2n}\cong
 \Gr_q^W \Ker \Big(H^{2n}_t\overset{N}{\lto} H^{2n}_t\Big)
 \]
 for $q<2n$. If we have $X^{[n]}=\emptyset$, 
 \eqref{sncgr} implies $\Gr_q^WH^{2n}=0$ for $q\leqslant n$, therefore 
 \[
 \Gr_q^WH^{2n}_t
 =\bigoplus _{i=0}^{\lfloor q/2\rfloor} 
 \Gr _{q-2i}^W \Ker \Big(H^{2n}_t\overset{N}{\lto} H^{2n}_t\Big)=0
 \] 
 for $q\leqslant n$, i.e., $N_{2n}^n=0$.
 By Theorem \ref{main}, we get $N_{2n}=0$, which is a contradiction.
\end{proof}

This corollary means in particular that there is no chain degeneration 
nor cycle degeneration of irreducible symplectic manifold, i.e. no good
degeneration such that the dual graph of the singular fiber is as
following: 

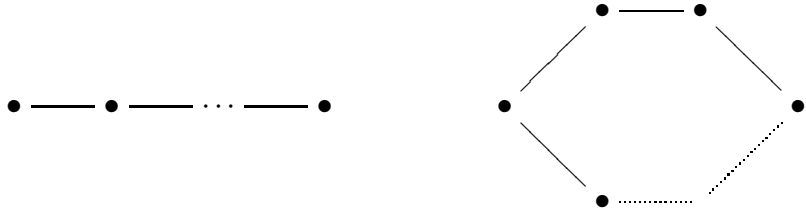
\[
\xymatrix{
&&&&&& \bullet \ar@{-}[r]  & \bullet \ar@{-}[rd]\\
\bullet \ar@{-}[r] & \bullet \ar@{-}[r]
 & \cdots \ar@{-}[r] &  \bullet 
 & & \bullet \ar@{-}[ru]\ar@{-}[rd] &  & & \bullet\\
&&&&&& \bullet \ar@{.}[r] & \ar@{.}[ru]
}
\]

\noindent
It is remarkable to compare this to the fact that 
we have cycle degenerations for generic degeneration of abelian
varieties and one can also generically
expect chain degenerations for Calabi-Yau manifolds.
We can regard the situation in the example of
Theorem \ref{exampleofdsv} is the least degenerate case of the
degeneration of symplectic manifolds. We should also note that 
we can expect some special property of the period map on the middle
cohomology of irreducible symplectic manifolds.


\end{document}